\documentclass{article}
\def\uu{\bigsqcup}

\def \Z {{\mathbf {Z}}}

\def \N {{\mathbf {N}}}
\def \J {{\mathbf {J}}}
\def \A {{\mathcal {A}}}
\def \B {{\mathcal {B}}}
\title{  Квазиподобие,  энтропия и дизъюнктность   эргодических систем}
\author{ В.В. Рыжиков, Ж.-П. Тувено}
\date{}

\usepackage[T1]{fontenc} % иногда T1 вместо T2A
\usepackage[cp1251]{inputenc}  
\usepackage[russian]{babel}
\begin{document}

\maketitle

\begin{abstract}   
В статье дан ответ  на вопрос А.М. Вершика о связи отношения  квазиподобия динамических систем с  энтропией Колмогорова. Доказано, что  все бернуллиевские действия заданной бесконечной счетной группы лежат в одном классе квазиподобия.  Открыт вопрос: попадает ли в этот  класс   небернуллиевское действие?   Антиподом  квазиподобия является  дизъюнктность (независимость) действий. М.С. Пинскер доказал, что детерминированное действие независимо от  действия с вполне положительной энтропией.   При помощи джойнингов в статье получено следующее обобщение теоремы Пинскера:  действие с  нулевой $P$-энтропией, инвариантом предложенным  А.А. Кирилловым и А.Г. Кушниренко, и действие с вполне положительной $P$-энтропией независимы.  
\end{abstract}
 
\Large
\section{Введение}

В эргодической теории известны несколько типов эквивалентности сохраняющих меру действий. К ним относятся изоморфизм, слабый изоморфизм, квазиподобие и спектральный изоморфизм.  
Под изоморфизмом подразумевается сопряжение  действий обратимым   преобразованием, сохраняющим меру. Слабый изоморфизм
двух систем означает, что каждая система изоморфна подсистеме, иначе говоря,  фактору другой системы. Напомним, что фактор -- это  ограничение действия на некоторую инвариантную сигма-алгебру. Квазиподобие двух действий означает
   существование инъективного марковского оператора с плотным образом, сплетающего эти действия. Отсутствие нетривиального марковского сплетения означает дизъюнктность систем.

Спектральный изоморфизм подразумевает наличие унитарного сплетающего оператора.
Все бернуллиевские автоморфизмы  спектрально изоморфны
(имеют счетнократный лебеговский спектр), но, как показал A.Н. Колмогоров,
их различает энтропия \cite{K}. Я.Г. Синай  доказал, что бернуллиевские автоморфизмы  с одинаковой энтропией    слабо изоморфны \cite{Si}, а Д. Орнстейн установил их изоморфизм \cite{O}.

В работе \cite{V} А.М. Вершик предложил понятие квазиподобия,
в связи с которым возник ряд интересных задач. Одна из них
    решена К. Франчеком и М. Леманчиком \cite{FL}: они предъявили квазиподобные автоморфизмы, не являющиеся слабо изоморфными. Вопрос об инвариантности энтропии относительно квазиподобия оставался открытым. Мы  отвечаем на него  отрицательно, показав, что все бернуллиевские действия заданной бесконечной счетной группы  квазиподобны. Доказательство  использует идею бернуллиевских джойнингов из работы \cite{T}.

Пусть действие $T$ обладает инариантной сигма-алгеброй, ограничение $T$ на эту алгебру называется фактором. Мы говорим $S$-фактор, если  фактор изоморфен 
действию $S$.
Свойство динамической системы называется \it наследственно устойчивым \rm, если 
оно наследуется как факторами, так  и  всякой системой, порожденной своими факторами с этим свойством. Инвариант "иметь нулевую энтропию" является наследственно устойчивым.
 М.С. Пинскер  в \cite{P} установил, что детерминированное действие дизъюнктно с К-системой. Это означает, что фактор с нулевой энтропией и фактор с вполне положительной энтропией всегда независимы. Это свойство равносильно  дизъюнкности
таких действий.  Мы предлагаем обобщение результата Пинскера, заменяя энтропию Колмогорова на  $P$-энтропию   Кириллова-Кушниренко (см. \cite{Ki}--\cite{R21}). 
С помощью  джойнингов  мы покажем, что   действие с вполне положительной 
$P$-энтропией  дизъюнктно с действием нулевой $P$-энтропии.    Примеры детерминированных $Z$-действий с вполне положительной $P$-энтропией будут 
предъявлены в классе   пуассоновских надстроек. 

\newpage
\section{Квазиподобие бернуллиевских действий}
%\vspace{3mm}

\bf Теорема 2.1. \it Все бернуллиевские действия  бесконечной счетной группы квазиподобны. \rm

\vspace{3mm}
Доказательство. Пусть $S$ — схема Бернулли типа 
$$\left(\frac 1 2, \, a,\, \frac {1} 2 - a\right), \ \ 0<a<1/2.$$
Она имеет следующие очевидные бернуллиевские факторы  $S_P$, $S_Q$, где $S_P$ имеет тип $\left(a, 1 -a\right)$, а $Q$ имеет тип $\left(\frac 1 2 , \frac 1 2\right)$.
Эти факторы порождают нашу систему $S$.

Идея доказательства следующая.  Пусть $E^Q$ — ортогональная проекция на пространство множителя $S_Q$
(это пространство обозначается $L_2(Q)$).
Определим оператор $\J$, установив $\J=E^Q\,|_{L_2(P)}$. Тогда $\J$ — марковский инъективный оператор, сплетающий факторы $S_P$ и $S_Q$, причем $\J L_2(P)$ плотен в $ L_2(Q)$. Итак, факторы $S_P$, $S_Q$ квазиподобны.
Построение такого оператора $\J$  подсказано  леммой 4 \cite{T}.

Теперь приведем  подробные рассуждения.
    Рассмотрим разбиения $$\xi=\{A, X\setminus A\},
\ \ 0<\mu(A)=a<1/2,$$ и $$\beta=\{ B, X\setminus B\},
\ \ \mu(B)=1/2, $$ где $$X=[0,1], \ \ A=[0,a], \ \ B=[0,1/2]. $$
Определим марковский оператор $J:L_2(X)\to L_2(X)$ как интегральный оператор со следующим $\xi\times\beta$-измеримым ядром $K(x,y)$:
$$ K(x,y)=0, \ \ (x,y)\in A\times B; \ \ \ \ \  K(x,y)=\frac 1{1-a},
\ \ (x,y)\in (X\setminus A)\times B;$$
$$ K(x,y)=2,
\ \ (x,y)\in A\times (X\setminus B); \ \ \ \ \  K(x,y)=\frac {1-2a}{1-a},
\ \ (x,y)\in (X\setminus A)\times (X\setminus B).$$
Напомним, что марковский оператор $J$ по определению сохраняет положительность функций, причем  $J$ и $J^\ast$ оставляют константы неподвижными.

Пусть $S$ — стандартный сдвиг в пространстве $X^\Z$, пусть $\A$ обозначает сигма-алгебру, порожденную множествами
$$ S^n(\dots\times X \times X \times A \times X \times X \dots),$$
и $\B$ -- сигма-алгебра, порожденная
$$ S^n(\dots\times X \times X \times B \times X \times X \dots).$$

Ограничение $S$ на
  $\A$ является автоморфизмом Бернулли $T_a$ с энтропией 
$$H(\xi)=-a\log_2 a -(1-a)\log_2(1-a),$$ а ограничение $S$ на
  $\B$ есть автоморфизм Бернулли $T$ с энтропией $H(\beta)=1$.

Теперь рассмотрим марковский оператор $\J_a: L_2(X^\Z,\mu^Z)\to L_2(X^\Z,\mu^Z)$
настройка $$ \J_a= \bigotimes_{z\in \Z} J =
  \dots\otimes J \otimes J \otimes J \otimes \dots.$$
Оператор $\J_a$ сплетает $T_a$ с автоморфизмом $T$:
$$ T\J_a=\J_a T_a.$$
Образ $\J_a$ плотен:
$$ \overline{\J_a L_2(\A)}= L_2(\B),$$
и
$$ Ker\J_a={0},$$
поскольку тензорные степени инъективного оператора $J$ также инъективны.
Таким образом, для всех $a$, $0<a<1/2,$ получаем, что  $T_a$ квазиподобен автоморфизму $T$.
 
Поскольку степень $T_a^n$ квазиподобена степени  $T^n_{a'}$, $0<a'<1/2,$
а для энтропии $h$ выполнено  $h(T_a^n)=nh(T_a)$, мы очевидным образом 
  для любых $c, c'>0$ найдем бернуллиевские автоморфизмы  $T$, $T'$ такие, что  $h(T)=c$, $h(T')=c'$ и $T$ квазиподобен  $T'$.

Замечаем, что оператор $$ \J_a\otimes \J_{a^2}\otimes \J_{a^3}\otimes\dots$$
сплетает автоморфизм с конечной энтропией с
автоморфизмом  $T\otimes T\otimes T\otimes\dots$, имеющим бесконечную энтропию. Применение теоремы Орнстейна об изоморфизме  \cite{O} завершает доказательство для $\Z$-действий.

Для  действий счетной группы $G$   повторяем приведенные выше рассуждения (подставляя $(X^G,\mu^G$) вместо $(X^\Z, \mu^\Z)$). Для  группы с  элементом бесконечного порядка используем результаты А.М. Степина \cite{St}. Для периодических счетных  бесконечных групп доказательство завершает применение  теоремы 1.1  Б. Сьюарда из \cite{Sw}.

\vspace{3mm}
Возникает интригующий вопрос: \it найдется  ли пара действий с вполне положительной энтропией, не являющихся квазиподобными?\rm

\section{Наследственно устойчивые  инварианты и $P$-энтропия}
Напомним, что инвариант динамической системы 
 называется  наследственно устойчивым, если он наследуется факторами и наследуется  системой, порожденной факторами, обладающими этим инвариантом.
  Заметим, что свойство иметь вполне положительную энтропию  не является наследственно устойчивым. В \cite{ST} доказано, что для заданного эргодического автоморфизма $T$ с нулевой энтропией найдется действие, порожденное двумя  бернуллиевскими факторами и при этом   обладающее $T$-фактором.

\vspace{3mm}
\bf Теорема 3.1. \it Пусть действие $S$ обладает наследственно устойчивым свойством,  а действие  $T$ таково, что оно и любой его нетривиальный фактор этим свойством не обладают. Тогда $S$ и $T$ дизъюнктны. \rm

\vspace{3mm}
Доказательство. Если действия  $S$ и $T$ не дизъюнктны,
существует система, порожденная счетным семейством $S$-факторов, имеющая нетривиальный фактор, изоморфный некоторому фактору действия $T$. 
Это утверждение в терминах джойнингов см. в \cite{T2}, лемма 3.2. Из условий теоремы мы видим, что такой нетривиальный фактор в нашем случае не   существует. Следовательно, действия $S$ и $T$ дизъюнктны.

\vspace{3mm}
\bf $P$-энтропия Кириллова-Кушниренко. \rm Рассмотрим  небольшую модификацию энтропии Кириллова-Кушниренко.    Пусть $P=\{P_j\}$ — последовательность конечных подмножеств в счетной бесконечной группе $G$. Для сохраняющего меру действия $T=\{T_g\}$ группы $G$ определим
$$h_j(T,\xi)=\frac 1 {|P_j|} H\left(\bigvee_{p\in P_j}T_p\xi\right),$$
$$h_{P}(T,\xi)={\limsup_j} \ h_j(T,\xi),$$
$$h_{P}(T)=\sup_\xi h_{P}(T,\xi),$$
%$$h^{inf}_{P}(T)=\sup_\xi\liminf_j \ h_j(T,\xi),$$
где $\xi$ обозначает конечное измеримое разбиение пространства $X$, а $H(\xi)$ — энтропия разбиения $\xi$:
$$ H(\{C_1,C_2,\dots, C_n\})=-\sum_{i=1}^n \mu( C_i)\ln \mu( C_i).$$

В статье Кириллова \cite{Ki} предполагается, что множества  $P_j$ образуют монотонную по включению последовательность.  Мы этого не требуем.  

Далее будет  рассмотрен  случай $G=Z$, причем для удобства множества  $P_j$ выбираются в виде  возрастающих прогрессий
    $$P_j=\{j,2j,\dots, L(j)j\}, \ \ \ L(j)\to\infty.$$

\bf Примеры автоморфизмов с нулевой $P$-энтропией. \rm
Пусть $T$ имеет нулевую энтропию. Степени $T^n$ также имеют нулевую энтропию.
Итак, для любых $\xi$ и $j$ существует $L(j)$ такой, что $$h_j(T^j,\xi)<\frac 1 j.$$
Мы можем выбрать такие $j_k\to \infty$, что для $P=\{P_{j_k}:k\in \N\}$ и всех конечных разбиений $\xi$ получим ${\limsup_k} \ h_{j_k }(T,\xi) =0.$ , $h_P(T)=0$ (см. \cite{R21}).

\vspace{3mm}
\bf Лемма 3.2. \it Нулевая $P$-энтропия наследственно стабильна. \rm

\vspace{3mm}
Лемма очевидным образом следует из того, что $\xi_1\vee\dots\vee \xi_n$
является порождающим разбиением для системы $F_1\vee\dots\vee F_n$, если $\xi_m$, $1\leq m\leq n$, являются порождающим разбиениями для соответствующих факторов $F_m$.

\vspace{3mm}
 Отметим, что такое наследственно устойчивое свойство как ''обладать дискретным спектром'', имеет энтропийную природу (см. \cite{Ku}, теорема 4).
Интересно узнать, \it исчерпываются ли все наследственно устойчивые инварианты комбинациями обсуждаемых энтропийных инвариантов? \rm 

\section{ Пуассоновские надстройки  с вполне положительной $P$-энтропией}
Говорим, что действие имеет вполне положительную $P$-энтропию, если любой нетривиальный фактор этого действия имеет положительную $P$-энтропию.
Из теоремы 3.1 и леммы 3.2 получаем следующее утверждение.

\vspace{3mm}
\bf Теорема 4.1. \it Пусть действие $T$ имеет 
вполне положительную $P$-энтропию, а  действие $S$ имеет нулевую     $P$-энтропию. 
Тогда $T$ и $S$ дизъюнктны ($S$-фактор и $T$-фактор всегда независимы). \rm

\rm
\vspace{3mm}
Примеры детерминированных систем с вполне положительной  $P$-энтропией удобно искать среди пуассоновских надстроек $T_\circ$ над бесконечными преобразованиями $T$ ранга один.

\bf Конструкции ранга один. \rm Пусть заданы параметры 
   $h_1=1$, $r_j\geq 2$ и наборы натуральных  чисел
$$ \bar s_j=(s_j(1), s_j(2),\dots,s_j(r_j)), \ s_j(i)\geq 0, \ j\in\N. $$
Фазовое пространство $X$ для преобразования $T$, определенного ниже,  представляет собой объединение башнен $$X_j=\uu_{i=0}^{h_j-1} T^iB_j,$$
где  $T^iB_j$ — непересекающиеся полуинтервалы, называемые этажами. 

На этапе $j$ преобразование  $T$ определено как обычный перенос интервалов (этажей), но на верхнем этаже 
$T^{h_j-1}B_j$ башни $X_j$ преобразование пока не определено.
 Башня $X_j$ разрезается на $r_j$ одинаковых узких подбашен $X_{j,i}$ (они называются колоннами),    и над каждой колонной  $X_{j,i}$  добавляется  $s_j(i)$ новых этажей. Преобразование $T$ по-прежнему определяется как подъем на этаж выше, но  самый последний надстроенный этаж над колонной с номером  $i$ преобразование $T$ отправляет в нижний этаж колонны   с номером  $i+1$.  Таким обазом возникает   новая башня 
$$X_{j+1}=\uu_{i=0}^{h_{j+1}-1} T^iB_{j+1}$$
высоты 
   $$h_{j+1}=r_jh_j+\sum_{i=1}^{r_j}s_j(i),$$
где $B_{j+1}$ -- нижний этаж колонны $X_{j,1}$.
Отметим, что доопределяя преобразование, мы полностью сохраняем предыдущие построениея. 
  Продолжая этот процесс до бесконечности,  получаем обратимое преобразование $T:X\to X$, сохраняющее меру  Лебега на объединении $X=\uu_j X_j$.
В \cite{R20},\cite{R23} можно найти разнообразные  применения таких конструкций, в том числе для гауссовских и пуассоновских надстроек.

\vspace{3mm}
\bf Пространство Пуассона. \rm Рассмотрим конфигурационное пространство $X_\circ$, состоящее из всех бесконечных счетных множеств $x_\circ$ таких, что каждый  интервал из пространств $X$ содержит лишь конечное число элементов множества $x_\circ$.

Пространство $X_\circ$ оснащается  мерой Пуассона. Напомним ее определение.
Подмножествам $A\subset X$ конечной $\mu$-меры в конфигурационном пространстве $X_\circ$  сопоставлены цилиндрические  множества
   $C(A,k)$, $k=0,1,2,\dots$,    по формуле
$$C(A,k)=\{x_\circ\in X_\circ \ : \ |x_\circ\cap A|=k\}.$$

Всевозможные конечные пересечения вида $\cap_{i=1}^N C(A_i,k_i)$
     образуют полукольцо.
На этом полукольце определена мера  $\mu_\circ$ следующим образом:
   при условии, что измеримые множества $A_1, A_2,\dots, A_N$ не пересекаются
и имеют конечную меру, положим
$$\mu_\circ(\bigcap_{i=1}^N C(A_i,k_i))=\prod_{i=1}^N \frac {\mu(A_i)^{k_i}}{k_i!} e ^{-\mu(A_i)}.\eqno (\circ)$$
Пояснение: если множества $A$, $B$ не пересекаются, то
    вероятность $\mu_\circ(C(A,k))\cap C(B,m))$  одновременного появления $k$ точек конфигурации $x_\circ$ в $A$
     и $m$ точек конфигурации $x_\circ$ в $B$ равна произведению вероятностей $\mu_\circ(C(A,k))$ и $\mu_\circ( C(B,m))$.
    Другими словами, события $C(A,k))$ и $C(B,m)$ независимы. Так как  множества $A_1, A_2,\dots,A_N$ не пересекаются,     в формуле $(\circ)$ фигурирует произведение соответствующих вероятностей. 
 Классическое продолжение меры является  пространством Пуассона
$(X_\circ,\mu_\circ)$,     изоморфным стандартному вероятностному пространству Лебега.
Автоморфизм $T$ пространства $(X,\mu)$ естественным образом индуцирует автоморфизм
$T_\circ$ пространства $(X_\circ,\mu_\circ)$,  называемый пуассоновской надстройкой над $T$.

Таким образом, c сохраняющим сигма-конечную меру действием группы 
ассоциирован ее   пуассоновский образ в группе автоморизмов вероятностного пространства. Этот замечательный факт нашел разнообразные  применения как  в эргодической теории, так и в теории представлений 
(см., например, \cite{VGG}-\cite{J}).

\vspace{2mm}
\bf Пуассоновские надстройки с вполне положительной $P$-энтропией. \rm
Пусть дана последовательность $L(j)\to\infty.$
Рассмотрим преобразования $T$  ранга один, для параметров которого выполнено  
$s_j(i)>L(j)h_j$.
Тогда для соответствующих башен $X_j$  выполнено $\mu (X_j)\to\infty$, причем  множества
$$X_j, \ T^{h_j}X_j,\ T^{2h_j}X_j,\ \dots,\ T^{L(j)h_j}X_j$$ не пересекаются.

Пусть $C=C(A,k)$, где $A\subset X_{j_0}$, для  пуассоновской надстройки $T_\circ$
замечаем, что множества
$$C, \ T^{h_j}_\circ C,\ T^{2h_j}_\circ C\ \dots, \ T^{L(j)h_j}_\circ C$$
$\mu_\circ$-независимы, так как $A, \ T^{h_j}A, \dots, T^{L(j)h_j}A$ не пересекаются.
Стандартные рассуждения показывают, что $T_\circ$ имеет вполне положительную $P$-энтропию  в случае 
$$P=\{ P_j \}, \ \ P_j=\{h_j, 2h_j,\dots, L(j)h_j\}.$$
Но при этом  $T_\circ$ как надстройка  над преобразованием ранга один  имеет нулевую энтропию (см. \cite{J}). 
Таким образом,  мы приходим к следующему утверждению.

\vspace{3mm}
\bf Теорема 4.2. \it  Для всякой детерминированной пуассоновской надстройки $S_\circ$  найдутся последовательность $P$, для которой $h_P(S_\circ)=0$, и  класс  детерминированных надстроек $T_\circ$ с вполне положительной  
$P$-энтропией.  Такие $S_\circ$ и $T_\circ$  дизъюнктны.   \rm

\vspace{3mm}
 Отметим, что аналогичная ситуация имеет место  для гауссовских действий.
Мы ограничились примерами  действий группы $\Z$, однако,  методы работы \cite{R21}  и настоящей статьи, несомненно,   позволяют строить  соответствующие примеры   для действий  счетных аменабельных групп. Подобный подход также может быть  полезным  для исследования масштабированной энтропии (см. \cite{VVZ}) детерминированных гауссовских и пуассоновских действий.

\vspace{2mm}
Авторы благодарят рецензента за замечания.

\large

\end{document}